\documentclass[final,narrowdisplay]{elsart1p}
\usepackage{amssymb}
\usepackage{amsmath}
\usepackage[notref,notcite]{showkeys}
\usepackage{graphicx}
\usepackage{psfrag}
\numberwithin{equation}{section}
\numberwithin{equation}{section}
\newtheorem{definition}{Definition}[section]

\newtheorem{Claim}[definition]{Claim}

\newtheorem{Remark}[definition]{Remark}
\newtheorem{Cor}[definition]{Corollary}

\newcommand \del    {\partial }
\newcommand \calR   {\mathcal R}
\newcommand \be     {\begin{equation}}
\newcommand \ee     {\end{equation}}
\newcommand \EQU        {\mathnormal{=}}

\newcommand \GT     {\mathnormal{>}}
\newcommand \LT     {\mathnormal{<}}
\newcommand \LEQ        {\mathnormal{\leq}}
\newcommand \GEQ        {\mathnormal{\geq}}
\newcommand \Se         {\mathbb{S}}
\newcommand \R      {\mathbb{R}}
\newcommand \Reee   {\mathbb{R}^3}
\newcommand \Rone   {\mathbb{R}^1}
\newcommand \See    {\mathbb{S}^2}
\newcommand \bfi    {\mathbf{i}}
\newcommand \bfn    {\mathbf{n}}
\newcommand \bfx    {\mathbf{x}}
\newcommand \bfF    {\mathbf{F}}

\newcommand \bfy    {\mathbf{y}}
\newcommand \bft    {\mathbf{t}}
\newcommand \bfPhi  {\mbox{\boldmath$\Phi$}}
\newcommand \bfnu   {\mbox{\boldmath$\nu$}}

\newcommand \scscL  {\!\scriptscriptstyle{L}}
\newcommand \scscR  {\!\scriptscriptstyle{R}}
\newcommand \scscT  {\!\scriptscriptstyle{T}}
\newcommand \uL         {{u}_{\scscL}}
\newcommand \uR         {{u}_{\scscR}}
\newcommand \uLlambda {{u}_{{\scscL},\lambda}}
\newcommand \uRlambda {{u}_{{\scscR},\lambda}}
\newcommand \uLphi  {{u}_{{\scscL},\phi}}
\newcommand \uRphi  {{u}_{{\scscR},\phi}}
\newcommand \diff   {{\mathrm{d}}}

\newcommand{\figref}[1]{Figure~\ref{#1}}

\newcommand{\EQSPLITLAM}
                     {\mbox{\eqref{eq:split-lambda.phi}$_{\lambda}\;$}}
\newcommand{\EQSPLITPHI}{\mbox{\eqref{eq:split-lambda.phi}$_{\phi}\;$}}

\begin{document}
\begin{frontmatter}

\title{Hyperbolic conservation laws on the sphere.
\\
A geometry-compatible finite volume scheme}

\vskip.9cm

\author{Matania Ben-Artzi, Joseph Falcovitz,}

\vskip.3cm

\address{Institute of Mathematics, Hebrew University,
\\
Jerusalem 91904, Israel.
\\
E-mail: {\tt mbartzi@math.huji.ac.il, ccjf@math.huji.ac.il}}

\vskip.15cm

\author{Philippe G. LeFloch}

\vskip.3cm

\address{Laboratoire Jacques-Louis Lions,
\\
Centre National de la Recherche Scientifique,
\\
Universit\'e de Paris 6, 4 place Jussieu,
\\ Ê     75252 Paris, France.
E-mail: \texttt{LeFloch@ann.jussieu.fr}
}

\begin{abstract}
We consider entropy solutions to the initial value problem
associated with scalar nonlinear hyperbolic conservation laws
posed on the two-dimensional sphere.
We propose a finite volume scheme which relies on a web-like mesh
made of segments of longitude and latitude lines.
The structure of the mesh  allows for a discrete version
of a natural geometric compatibility condition, which arose earlier
in the well-posedness theory established by Ben-Artzi and LeFloch.
We study here several classes of
flux vectors which define the conservation law under consideration.
They are based on prescribing a suitable vector field in the
Euclidean three-dimensional space and then
suitably projecting it on the sphere's tangent plane;
even when the flux vector in the
ambient space is constant, the corresponding flux vector is a
non-trivial vector field on the sphere.
In particular, we construct here ``equatorial periodic solutions'',
analogous to one-dimensional periodic solutions to one-dimensional
conservation laws,
as well as a  wide variety of stationary (steady state) solutions.
We also construct ``confined solutions'', which are time-dependent
solutions supported in an arbitrarily specified subdomain of the sphere.
Finally, representative numerical examples and test-cases are presented.
\end{abstract}

\end{frontmatter}

\tableofcontents

\section{Introduction}
\label{IN-0}

In this paper, building on our earlier analysis in \cite{BL,ABL}
we study in detail
the class of scalar hyperbolic conservation laws posed on the
two-dimensional unit sphere
$$
\See=\left\{(x,y,z)\in\Reee,\; x^{2}+y^{2}+z^{2}=1\right\}.
$$
We propose a Godunov-type finite volume scheme that satisfies
certain important consistency and convergence properties.
We then present a
second-order extension based on the generalized Riemann problem
(GRP) methodology~\cite{grp_book}.

It should be
stated at the outset that an important motivation for this paper
is the need to
provide accurate numerical tools for the so-called shallow water
system on the sphere.
This system is widely  used in geophysics as a model for global air
flows on the rotating Earth~\cite{Haltiner-1971}.
In its mathematical classification it is a
system of nonlinear hyperbolic PDE's posed on the sphere.
Its physical nature dictates that it can be
described ``invariantly'', namely in a way which is independent of
any particular
coordinate system. Locally, it has the (mathematical) character of a
two-dimensional
isentropic compressible flow, whereas globally the spherical
geometry plays a crucial
role in shaping the nature of solutions --which, as expected
for nonlinear
hyperbolic equations, may contain propagating discontinuities such
as shock fronts or contact curves.
Thus, the relation of the present study to the shallow water system  is
analogous to the connection between Burgers' equation and the
system of compressible fluid flow (say, in the plane).
In fact, in light of this analogy it is somewhat
surprising that in the existing literature so far, virtually all
treatments, theoretical as well as numerical, were confined to the Cartesian setting.
In particular, to the best of our knowledge,
there have been no systematic numerical studies of scalar conservation laws
on the sphere.

Having introduced the scalar conservation law as a simple model for
more complex physical
systems, we should emphasize here also the intrinsic mathematical
interest of the model under consideration.
It is already known (see~\cite{BFLi} and references there) that even in
the Cartesian setting, the two-dimensional scalar conservation law
displays a wealth
of wave interactions typical of the physical phenomena
(such as triple points, sonic
shocks, interplay of rarefactions and shocks coming from different
directions and more).
As we show here, ``geometric effects'', superposed on the (necessarily)
two-dimensional framework, carry the scalar model still further.
For example, the
concept of ``self-similar'' solutions makes no sense here.
In particular, one loses the
Riemann solutions, a fundamental building block in many schemes
(of the so-called
``Godunov-type"). On the other hand, it allows for large classes of
non-trivial steady
states, periodic solutions, and solutions supported in specified
subdomains.
All these
have natural consequences in developing numerical schemes;
they offer us a variety of
test-cases amenable to detailed analysis, to be compared with the
computational results.

In practical applications a finite volume scheme requires a
specification of a coordinate system,
where the symmetry-preserving latitude--longitude coordinates are
the ``natural coordinates'' of preferred choice.
The proposed finite volume scheme in this paper
is based on these natural
coordinates, but should pay attention to the artificial singularities
at the poles.

In \cite{ABL}, a general convergence theorem was proved for a class
of finite volume
schemes for the computation of entropy solutions to conservation
laws posed
on a manifold. As a particular example, the case of the sphere
$\See$ was
discussed, both from the points of view of an ``invariant'' formalism
and that of an ``embedded'' coordinate-dependent formulation.
In the present study we focus on the
sphere $\See$ and we actually construct, in a fully explicit and
implementable way, a finite volume scheme which is geometrically
natural and can be
viewed as an extension of the basic Godunov scheme for
one-dimensional conservation laws.
Furthermore, we prove that our scheme fulfills all of the assumptions required in
\cite{ABL}, which ensures its strong convergence toward the unique
entropy solution to the initial value problem under consideration.
We then describe the GRP extension of
the scheme, whose convergence proof is still a challenging open problem.

The theoretical background about the well-posedness theory for
hyperbolic conservation laws on manifolds
was established recently by Ben-Artzi and LeFloch~\cite{BL}
together with collaborators \cite{ALO,ABL,LNO}.
An important condition arising in the theory is
the ``zero-divergence'' or {\sl geometric-compatibility} property of
the flux vector;
a basic requirement in our construction of a finite volume scheme
is to formulate and ensure a suitable discrete version of this condition.

\begin{figure}[!htb]
\psfrag{Lam}[][ml]{\hspace{ 1.3em}{$\lambda$}}
\psfrag{Phi}[][ml]{\hspace{ 0.8em} {$\phi$}}
\begin{center}
\includegraphics[clip,width =80mm]{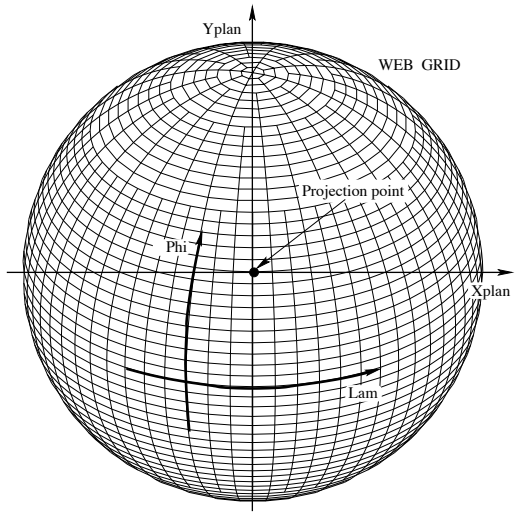}
\caption{Web grid on a sphere}
\label{fig:sphere.1}
\end{center}
\end{figure}

We conclude this introduction with some notation and remarks
connecting the present
paper to the general finite volume framework presented in~\cite{ABL}.
Following the terminology therein, we use an ``embedded'' approach
to the spherical geometry,
namely, we view the sphere as embedded in the three-dimensional
Euclidean space $\Reee$.
We denote by $\bfx$ a variable point on the sphere $\See$, which can be
represented in terms of its longitude $\lambda$ and its latitude $\phi$.
Following the
conventional notation in the geophysical literature we assume that
$$
0 \leq\lambda\leq 2\pi, \quad -\frac{\pi}{2} \leq \phi\leq\frac{\pi}{2},
$$
so that the ``North pole'' (resp. ``South pole'') is at
 $\phi=\frac{\pi}{2}$ (resp.  $-\phi=\frac{\pi}{2}$) and the equator
is $\big\{ \phi=0, \, 0\leq\lambda\leq 2\pi \big\}.$
(See \figref{fig:sphere.1}.)
The coordinates in $\Reee$ are denoted by
$(x_{1},x_{2},x_{3}) \in \Reee$ and the corresponding unit vectors are
$\bfi_1,\,\bfi_2,\,\bfi_3$.
Thus, at each point $\bfx=(\lambda,\phi)\in\See$, the unit tangent
vectors (in the $\lambda,\phi$ directions) are given by
\be\begin{split}
\bfi_\lambda &= -\sin\lambda         \,\bfi_1
                 +\cos\lambda         \,\bfi_2,\\
\bfi_{\phi   }&=  \hspace{0pt}
                 -\sin\phi\,\cos\lambda\,\bfi_1
                 -\sin\phi\,\sin\lambda\,\bfi_2
                 +\cos\phi            \,\bfi_3.
\label{eq:sphere.2}
\end{split}
\nonumber
\ee
It should be observed that while a choice of a coordinate system is
necessary in practice, it always introduces {\sl singularities} and
the unit vectors given above are {\sl not} well-defined at the poles
and, therefore, in the neighborhood of these points it
cannot be used for a representation of
smooth vector fields (such as the flux vectors of our conservation laws).
We also emphasize that the status of these two poles is equivalent
to the one of any other pair of opposite points on the sphere.
When such local coordinates are introduced, special care is needed
to handle these points in practice, and this is precisely why we advocate
a different approach.

Continuing with the description of our ``embedded'' approach, we define
the unit normal  $\bfn_{\bfx}$, to $\See$ at some point $\bfx$ by
\be
\begin{split}
\bfn_{\bfx   }&=  \cos\phi\,\cos\lambda\,\bfi_1
                 +\cos\phi\,\sin\lambda\,\bfi_2
                 +\sin\phi            \,\bfi_3.
\label{eq:sphere.3}
\end{split}
\nonumber \ee
Then, any tangent vector field $\bfF$ to $\See$ is represented by \be
\begin{split}
\bfF&=  F_\lambda \, \bfi_\lambda + F_\phi \, \bfi_\phi
\label{eq:sphere.4}
\end{split}
\nonumber
\ee
and the tangential gradient operator is
$$ 
\nabla_T= \left( \frac{1}{\cos\phi}\frac{\del}{\del
                   \lambda},\frac{\del}{\del\phi} \right).
$$
Thus, the (tangential) gradient of a scalar function
$h(\lambda,\phi)$ is given by
\be
\label{eq:sphere.4.2}
\nabla_T h=\frac{1}{\cos\phi}\frac{\del h}{\del
           \lambda} \, \bfi_\lambda+\frac{\del h}{\del\phi} \,
            \bfi_\phi,
\ee
and the divergence of a vector field $\bfF$ is
\be
\begin{split}
\nabla_{\scscT}\cdot\bfF &=
  \frac{1}{\cos\phi}\left(
  \frac{\del}{\del\phi}\left(F_\phi\,\cos\phi\right)
 +\frac{\del}{\del\lambda}F_\lambda\right).
\label{eq:sphere.5}
\end{split} 
\ee

Given now a vector field $\bfF = \bfF(\bfx, u)$ depending on a real
parameter $u$, the
associated hyperbolic conservation law under consideration is
\be
\begin{split}
&\frac{\del u}{\del t} + \nabla_{\scscT}\cdot
  \Big( \bfF(\bfx,u)\Big) = 0,
  \qquad (\bfx,t)\in\See\times[0,\infty),
\label{eq:sphere.6}
\end{split}
\ee where $u=u(\bfx,t)$ is a scalar unknown function, subject to
the initial condition
\be
\begin{split}
& u(\bfx,0) = u_0(\bfx), \qquad \bfx\in\See
 \label{eq:sphere.7}
\end{split}
\ee
for some prescribed data $u_0$ on the sphere.
 As mentioned above, we will impose on the vector field
$\bfF(\bfx, u)$ an additional ``geometry compatibility'' condition.

An outline of this paper is as follows.
In Section~\ref{constr}, we consider the construction of geometry-compatible flux vectors, while
Section~\ref{classf} is devoted to a description
of several families of special solutions associated with the
constructed flux vectors.
In Section~\ref{Description.sec}
we discuss our (first-order) finite volume scheme, which can be
regarded as
a Godunov-type scheme. We prove that it satisfies all of the
assumptions imposed on
general finite volume schemes in~\cite{ABL}, and we conclude that
it converges to the
exact (entropy) solution. In Section~\ref{second} we describe the
(second-order) GRP extension of the scheme.
Finally, in Section~\ref{numeri} we present a variety of numerical
test cases.

\section{Families of geometry-compatible flux vectors}
\label{constr}

As pointed out in \cite{ABL}, every smooth vector field
$\bfF(\bfx,u)$ on $\See$ can be represented in the form
\be
\begin{split}
&\bfF(\bfx,u)=\bfn(\bfx)\times\bfPhi(\bfx,u),
\label{eq:sphere.8}
\end{split}
\ee
where $\bfPhi(\bfx,u)$ is a restriction to $\See$ of a vector field
(in $\Reee$) defined in some neighborhood (i.e., a ``spherical shell'')
of $\See$ and for all values of the parameter $u$.
The basic requirement imposed now on the flux vector $\bfF(\bfx,u)$ is  the
following {\sl divergence free} or {\sl geometric compatibility}
condition:
For any fixed value of the parameter $v\in\R$,
\be
\nabla_{\scscT}\cdot\bfF(\bfx,v)= 0.
\label{eq:sphere.9}
\ee
A flux vector $\bfF(\bfx,u)$ satisfying \eqref{eq:sphere.9} is called
a {\sl geometry-compatible flux} \cite{BL}.
Note that this condition is equivalent, in
terms of the nonlinear
conservation law~\eqref{eq:sphere.6}, to the following
requirement: {\sl constant initial data are (trivial) solutions to
the conservation law.}
In the case of the sphere $\See$ the condition~\eqref{eq:sphere.9}
can be recast in
terms of a condition on the vector field $\bfPhi(\bfx,u)$ appearing
in~\eqref{eq:sphere.8}.
See \cite[Proposition 3.3]{ABL}.

Our main aim in the present section is singling out two (quite general)
families of geometry-compatible fluxes of particular interest, which are amenable to detailed
analytical and numerical investigation.

The flux-vectors of interest are introduced by way of the following two claims.

\

\begin{Claim}[Homogeneous flux vectors.]
\label{cor2.1}
If the three-dimensional flux $\bfPhi(\bfx,u)=\bfPhi(u)$ is
independent of $\bfx$ (in a neighborhood of $\See$), then the
corresponding flux vector
$\bfF(\bfx,u)$ given by~\eqref{eq:sphere.8} is geometry-compatible.
\end{Claim}

\

\noindent{\sl Proof.}
The following decomposition applies to any vector $\bfPhi(u) \in \Reee$ in the form
\be
\begin{split}
\bfPhi(u) & = f_1(u)\,\bfi_1 + f_2(u)\,\bfi_2 + f_3(u)\,\bfi_3,
\label{eq:sphere.10}
\end{split}
\ee
so that $\bfF(\bfx,u)=F_\lambda(\lambda,\phi,u) \, \bfi_\lambda+
              F_\phi(\lambda,\phi,u) \, \bfi_\phi$, with
\be
\begin{split}
F_\lambda(\lambda,\phi,u)&=  \hspace{9pt}
                               f_1(u)\,\sin{\phi}\,\cos{\lambda}
                              +f_2(u)\,\sin{\phi}\,\sin{\lambda}
                              -f_3(u)\,\cos{\phi},\\
F_{\phi   }(\lambda,\phi,u)&= -f_1(u)\,           \sin{\lambda}
               \hspace{24pt}  +f_2(u)\,           \cos{\lambda}.
\label{eq:sphere.11}
\end{split}
\ee
We can directly apply the divergence operator
\eqref{eq:sphere.5} to $\bfF(\bfx,u)$ and the desired claim follows.
\qed

\

\begin{Claim}[Gradient flux vectors.]
\label{cor2.2}
Let  $h=h(\bfx,u)$ be a smooth function  of the variables $\bfx$
(in a neighborhood of $\See$)
and $u\in\R$, and consider the associated three-dimensional
flux $\bfPhi(\bfx,u)=\nabla h(\bfx,u)$ (restricted to $\bfx\in\See$).
Then, the flux vector $\bfF(\bfx,u)$ given
by \eqref{eq:sphere.8} is geometry-compatible.
\end{Claim}

\

\noindent{\sl Proof.}
We use the divergence theorem in an arbitrary domain $D\subseteq\See$
 with smooth boundary $\del D$:
$$
\aligned
\int_D \nabla_{\scscT}\cdot \big( \bfF(\bfx,v) \big) \, d\sigma
& = \int_{\del D}\bfF(\bfx,v)\cdot \bfnu(\bfx) \, ds
\\
& = \int_{\del D} \big( \bfn(\bfx)\times\nabla h(\bfx,v) \big)
\cdot \bfnu(\bfx) \, ds,
\endaligned
$$
where $\bfnu(\bfx)$ is the unit normal (at $\bfx$) along
$\del D \subset \See$, $d\sigma$ is the
surface measure on $\See$, and $ds$ is the arc length along $\del D.$

In particular, $\bfn(\bfx)\times \bfnu(\bfx)=\bft(\bfx)$ coincides
with the (unit) tangent vector to $\del D$ at $\bfx$.
It follows that the triple product $\big( \bfn(\bfx)\times\nabla
h(\bfx,u) \big) \cdot \bfnu(\bfx)=\nabla h(\bfx,u) \cdot \bft(\bfx)$
is nothing but the directional
derivative $\nabla_{\del D}$ of $h$ along $\del D$. Since
$$
\int_{\del D} \nabla_{\del D} h \, ds = 0,
$$
we thus find
$$
\int_D \nabla_{\scscT}\cdot\bfF(\bfx,u) \, d\sigma=0,
$$
and since this holds for any smooth domain $D$, we conclude that
$\nabla_{\scscT}\cdot\bfF(\bfx,v)=0$ for all $v\in\R$.
\qed

\

\begin{Remark}
1. Claim~\ref{cor2.1} is a special case of Claim~\ref{cor2.2}.
Indeed, by
taking in the latter $h(\bfx,u) = x_1f_1(u)+x_2f_2(u)+x_3f_3(u)$
we obtain the conclusion
of the former. However, we chose to single out Claim~\ref{cor2.1}
as a special case since
it will serve in obtaining special solutions (Section~3) and in
dealing with numerical examples (Section~6).

2. The steps in the construction of the gradient flux vector in
Claim~\ref{cor2.2} are
``linear in nature'', namely if $h(\bfx,u) = h_1(\bfx,u) + h_2(\bfx,u)$
then the corresponding (geometry-compatible) flux vectors satisfy
$\bfF(\bfx,u)=\bfF_1(\bfx,u)+\bfF_2(\bfx,u).$
However, it is clear that the
corresponding solutions to~\eqref{eq:sphere.6} do not add up linearly,
due to the nonlinear dependence in $u$.
\end{Remark}

\section{Special solutions of interest}
\label{classf}

\subsection{Periodic equatorial solutions}

The scalar conservation laws discussed in this paper have two basic
features:
\begin{itemize}
\item The problem is necessarily two-dimensional (in spatial
coordinates).
\item The geometry plays a significant role, inasmuch as the flux
vectors are subject to geometric constraints.
\end{itemize}
It should be noted that even within the framework of
\textit{Euclidean} two
dimensional conservation laws there is a great wealth of special
solutions,
displaying complex wave interactions, such as triple points,
sonic shocks and more.
We refer to~\cite{li_book,BFLi} for detailed treatments of the
theoretical and numerical aspects.

In the situation under consideration in the present paper, geometric
effects yield a large variety
of non-trivial steady states, solutions supported in arbitrary
subdomains, etc.
In this section we consider such solutions by selecting some special
flux vectors $\bfF(\bfx,u)$ on $\See$.
This is accomplished by making special choices of
$\bfPhi(\bfx,u)$ in the general representation (see~\eqref{eq:sphere.8})
$\bfF(\bfx,u)=\bfn(\bfx)\times\bfPhi(\bfx,u)$,
where $\bfPhi(\bfx,u)$ is a restriction to $\See$ of a vector field
(in $\Reee$)
defined in some neighborhood (i.e., ``spherical shell'') of
$\See$ and for all values of the parameter $u$.

We begin our discussion with the case of {\sl periodic equatorial
solutions,} defined as follows.
Taking $f_1(u)=f_2(u)\equiv 0$ in the general decomposition \eqref{eq:sphere.10} so that,
by~\eqref{eq:sphere.11},
\begin{equation*}
\begin{split}
F_\lambda(\lambda,\phi,u)&=
                              -f_3(u)\,\cos{\phi},\\
F_{\phi   }(\lambda,\phi,u)&= 0,
\end{split}
\end{equation*}
the conservation law~\eqref{eq:sphere.6} takes the particularly simple form
\be
\begin{split}
&\frac{\del u}{\del t} - \frac{\del}{\del\lambda}f_3(u) = 0,
                    \qquad (\bfx,t)\in\See\times[0,\infty).
\label{equator}
\end{split}
\ee
In particular, obtain the following important conclusion.

\

\begin{Cor}[Solutions with one-dimensional structure.]
Let $\widetilde u= \widetilde u(\lambda,t)$ be a solution to the
following one-dimensional conservation law with periodic boundary condition
$$
\frac{\del  \widetilde u}{\del t} - \frac{\del}{\del\lambda}f_3(
\widetilde u) = 0,
 \quad 0<\lambda\leq 2\pi,\quad
\widetilde u(0,t)= \widetilde u(2\pi,t),
$$
and let $\widehat u = \widehat u(\phi)$ be an arbitrary function.
Then, the function $u(\lambda,\phi, t)=
\widetilde  u(\lambda,t)\, \widehat u(\phi)$ is a solution to the
conservation law~\eqref{equator}.
\end{Cor}

\

It follows that all periodic solutions from the one-dimensional case can be recovered
here as special cases. However, in numerical experiments the
computational grid is two-dimensional, so it is not obvious
that the accuracy achieved in the computation of the former
can indeed be achieved in the numerical scheme implemented on the sphere.
This issue will be further discussed below, in Section~6.


\subsection{Steady states}

Let $\bfF = \bfF(\bfx,u)$ be a flux vector and $u_0: \See \to \R$ be
an initial function such that
${\nabla_{\scscT}\cdot \big( \bfF(\bfx,u_0(\bfx)) \big) \equiv 0}$.
Then, clearly $u_0$ is a
stationary solution (or steady state) to the conservation law.
In fact, we can show that there exist many (analytically computable)
non-trivial steady state solutions, as follows.

\

\begin{Claim}[A family of steady-state solutions.]
\label{steady-claim}
Let $h=h(\bfx,u)$ be a smooth function defined for all
$\bfx$ in a neighborhood of $\See$, and consider the associated gradient flux vector
$\Phi = \nabla h$ (as in Claim~\ref{cor2.2}).
Suppose the function $u_0: \See \to \R$ satisfies the condition
\begin{equation}
\nabla_{\bfy}h(\bfy,u_0(\bfx))|_{\bfy=\bfx}=
\nabla_{\bfx}H(\bfx),\qquad \bfx\in\See,
\end{equation}
where $H=H(\bfx)$ be a smooth function defined in a neighborhood of $\See$.
Then, $u_0$ is a stationary solution to the conservation law~\eqref{eq:sphere.6}.
\end{Claim}

\

\noindent{\sl Proof.}
We follow the proof of Claim~\ref{cor2.2} and the notation therein.
Using the divergence
theorem in an arbitrary domain $D\subseteq\See$ with smooth boundary
$\del D$, we obtain
\begin{equation*}
\aligned
\int_D \nabla_{\scscT}\cdot \big( \bfF(\bfx,u_0(\bfx)) \big) \, d\sigma
& = \int \limits_{\del D}\bfF(\bfx,u_0(\bfx))\cdot \bfnu \, ds
\\
& = \int \limits_{\del D} \big( \bfn(\bfx)\times\nabla_{\bfx}
  H(\bfx) \big) \cdot \bfnu(\bfx) \, ds.
\endaligned
\end{equation*}
where, as before, $\bfnu(\bfx)$ is the unit normal, $d\sigma$ the surface measure,
and $ds$ the arc length. In particular,
 $\bfn(\bfx)\times \bfnu(\bfx)=\bft(\bfx),$
the (unit) tangent vector to $\del D$ at $\bfx.$
It follows  that the triple product $(\bfn(\bfx)\times\nabla_{\bfx}
H(\bfx))\cdot \bfnu(\bfx)=(\nabla _{\bfx}H(\bfx))\cdot \bft(\bfx)$ is
the directional derivative of $H$ along $\del D.$
Thus,
$$
\int_D \nabla_{\scscT}\cdot \big( \bfF(\bfx,u_0(\bfx)) \big) \, d \sigma  =0,
$$
and since this holds for any smooth domain $D$, it follows that
$\nabla_{\scscT}\cdot \big( \bfF(\bfx,u_0(\bfx)) \big) \equiv 0,$ which concludes
the proof.
\qed

\

The above claim yields readily a large family of non-trivial stationary
solutions, as expressed in the following corollary.

\

\begin{Cor}\label{1Dsteady} Consider the flux vector
$\bfF=\bfF(\bfx,u)$ given by
$$
\bfF(\bfx,u)=\bfn(\bfx)\times \big( f_1(u)\,\bfi_1\big),
$$
for an arbitrary choice of function $f_1=f_1(u)$.
Then, any function $u_0=u_0(x_1)$ depending only on the first
coordinate $x_1$ is a stationary
solution to the conservation law (associated with this flux).
In particular, in polar coordinates $(\lambda,\phi)$ any function
of the form $u_0(\lambda,\phi)=g(\cos\phi\cos\lambda)$ is a
stationary solution.
\end{Cor}

\

\noindent{\sl Proof.}
According to Claim~\ref{cor2.1} this flux vector is associated with
the scalar function $h(\bfx,u)=x_1f_1(u).$
So we can invoke Claim~\ref{steady-claim} with
$H(\bfx)=H(x_1)$ such that $H'(x_1)=f_1(u_0(x_1)).$
\qed

\

\begin{Remark}
\label{band1}
This corollary enables us to construct stationary solutions supported
in ``bands'' on the sphere.
This is accomplished by taking $u_0=u_0(x_1)$ to be supported in
$0<\alpha<x_1<\beta<1.$
Observe that this band is not parallel neither to the latitude curves
($\phi=const$) nor to the longitude curves ($\lambda=const$).
\end{Remark}

\

There is yet another possibility of obtaining stationary solutions,
where all three coordinates are involved, as stated now.
This example can also be derived from the previous one by applying a rotation in $\R^3$.

\

\begin{Cor} Consider the flux vector $\bfF=\bfF(\bfx,u)$ be given by
$$
\aligned
\bfF(\bfx,u)
& = \bfn(\bfx)\times (f_1(u)\,\bfi_1+f_2(u)\,\bfi_2+f_3(u)\,\bfi_3)
\\
& = f(u) \, \bfn(\bfx)\times ( \bfi_1 + \bfi_2 + \bfi_3),
\endaligned
$$
in which all three components coincide: $f_1(u)=f_2(u)=f_3(u)=f(u)$.
Then, any function of the form
$u_0(\bfx) = \widetilde u_0(x_1+x_2+x_3)$,  where
$\widetilde u_0$ depends on one real variable, only,
is a stationary solution to the conservation law associated with
the above flux.
\end{Cor}

\

\noindent{\sl Proof.}
Following the proof of the previous corollary, we now take
$H(\bfx)=H_0(x_1+x_2+x_3),$ where $H_0'(\xi)=f(\widetilde u_0(\xi))$.
\qed

\

\begin{Remark}
\label{band2}
In analogy with Remark~\ref{band1}, this result allows us to construct
stationary solutions in a spherical ``cap'' (a piece of the sphere
cut out by a plane).
In Section~6 below, we will provide numerical test cases for such stationary solutions.
\end{Remark}

\subsection{Confined solutions}
If in the  conservation law~\eqref{eq:sphere.6} we have
$\bfF(\bfx,u)\equiv 0$ for
$\bfx$ in the \textit{exterior} of some domain $D\subseteq \See,$
identically in
$u\in\R,$ and if the initial function $u_0(\bfx)$ vanishes outside of
$D$, then
clearly the solutions satisfy $u(\bfx,t)=0$ for $\bfx\notin D$ and
all $t\geq 0.$ We label such
solutions as \textsl{confined (to $D$) solutions}. In view of
equation~\eqref{eq:sphere.8} a sufficient condition for the vanishing of
$\bfF(\bfx,u)$ outside of
$D$ is obtained by $\bfPhi(\bfx,u)=0$ for $\bfx\notin D,$
identically in $u\in\R.$
In view of Claim~\ref{cor2.2}, this will follow if we choose
$h(\bfx,u)$ such
that $h(\bfx,u)\neq 0$ for $\bfx$ only in $D.$ In particular, let $\psi=\psi(\xi)$ be
a twice continuously differentiable function on $\R$ supported in
the interval $(\alpha,\beta) \subseteq
     (0,1)$ and such
     that $3\beta^2>1$ and $3\alpha^2<1.$
      With an
     eye to computable test cases, we can use this function to generate
     solutions which are confined within
     the intersection of $\See$ with the (three-dimensional) cube
     $[\alpha,\beta]^3$.

\
\begin{Claim}[A family of confined solutions.]
Let $\psi$ be as above and let $f=f(u)$ be any (smooth) function of
$u\in\R$.  Define  $h=h(\bfx,u)$  by
$$ h(\bfx,u)=\psi(x_1) \psi(x_2) \, \psi(x_3) \, f(u),$$
and let  $\bfF(\bfx,u)$ be the gradient flux vector determined in terms
of $h(\bfx,u)$
as in Claim~\ref{cor2.2}.
Let $D\subseteq\See$ be the spherical patch cut out from $\See$
by the inequalities $\alpha<x_i<\beta,\quad i=1,2,3.$
Then, if the initial data
$u_0(\bfx)$ is supported in $D,$ the solution $u=u(\bfx,t)$ of the
conservation law~\eqref{eq:sphere.6}
associated with $\bfF(\bfx,u)$ is supported in $D$ for all $t\geq 0.$
\end{Claim}

\

Possible choices for a function $\psi: [\alpha, \beta] \to \R$ as in the claim are
$\psi(\xi)=\sin^2(k\xi)$ for some integer $k$ such that $k\alpha$
and $k\beta$
are multiples of $\pi$,
 or else $\psi(\xi)=(\xi-\alpha)^2(\xi-\beta)^2$.

\section{Design of the scheme}
\label{Description.sec}

\subsection{Computational grid}
\label{computational_grid.sec}

\begin{figure}[!htb]
\psfrag{Lambda}[][ml]{\hspace{ 2.0em} {$\lambda$}}
\psfrag{Lambda1}[][bl]{\hspace{ 2.5em} {$\lambda_1$}}
\psfrag{Lambda2}[][bl]{\hspace{ 2.5em} {$\lambda_2$}}
\psfrag{Lambdam}[][bl]{\hspace{ 2.5em}
       {$(\lambda^{e,m},\phi^{e,m})$}}
\psfrag{Phi}[][ml]{\hspace{ 0.7em} {$\phi$}}
\psfrag{Phi1}[][ml]{\hspace{ 1.5em} {$\phi_1$}}
\psfrag{Phi2}[][ml]{\hspace{ 1.5em} {$\phi_2$}}
\psfrag{Phim}[][ml]{\hspace{-3.2em}
       {$(\lambda^{e',m},\phi^{e',m})$}}
\psfrag{nu}[][ml]{\hspace{ 0.0em} {$\bfnu$}}
\psfrag{Rectangle}[][ml]{\hspace{ 2.0em} {$\calR$}}
\psfrag{DR}[][ml]{\hspace{ 1.0em} {$\del\calR$}}
\psfrag{E}[][ml]{\hspace{ 0.3em} {$e$}}
\psfrag{Eprime}[][tl]{\hspace{ 0.3em} {$e'$}}
\begin{center}
\includegraphics[clip,bb=57 266 540 572,width =80mm]
                                                    {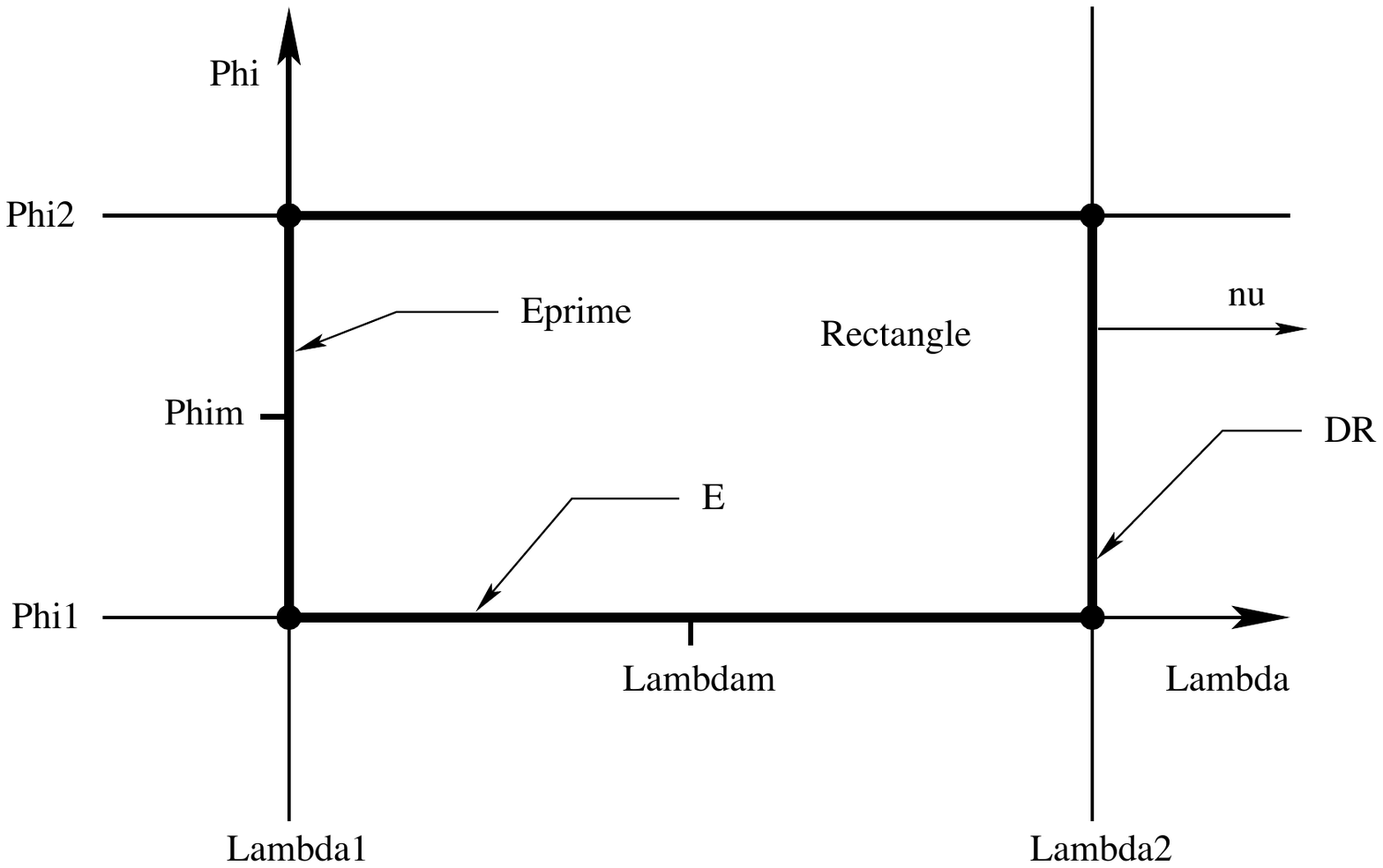}
\caption{Rectangular cell $\calR$ as part of grid on $\See$}
\label{fig:cells.1}
\end{center}
\end{figure}

\begin{figure}[!htb]
\psfrag{Lambda}[][ml]{\hspace{ 2.0em} {$\lambda$}}
\psfrag{Lambda1}[][bl]{\hspace{ 2.0em} {$\lambda_1$}}
\psfrag{Lambda2}[][bl]{\hspace{ 2.0em} {$\lambda_2$}}
\psfrag{Lambda3}[][ml]{\hspace{ 2.0em} {$\lambda_3$}}
\psfrag{Lambdam}[][bl]{\hspace{ 2.5em}
       {$(\lambda^{e,m},\phi^{e,m})$}}
\psfrag{Phi}[][ml]{\hspace{ 1.1em} {$\phi$}}
\psfrag{Phi1}[][ml]{\hspace{ 1.2em} {$\phi_1$}}
\psfrag{Phi2}[][ml]{\hspace{ 1.2em} {$\phi_2$}}
\psfrag{Phim}[][ml]{\hspace{ 1.5em}  {$\phi_{m}$}}
\psfrag{Phim}[][ml]{\hspace{-3.2em}
       {$(\lambda^{e',m},\phi^{e',m})$}}
\psfrag{nu}[][ml]{\hspace{ 0.0em}  {$\bfnu$}}
\psfrag{Rectangle}[][ml]{\hspace{ 2.0em}  {$\calR$}}
\psfrag{DR}[][ml]{\hspace{ 1.3em}  {$\del\calR$}}
\psfrag{E}[][ml]{\hspace{ 0.3em} {$e$}}
\psfrag{Eprime}[][tl]{\hspace{ 0.3em} {$e'$}}
\begin{center}
\includegraphics[clip,bb=57 191 540 654,width =80mm]
                                                  {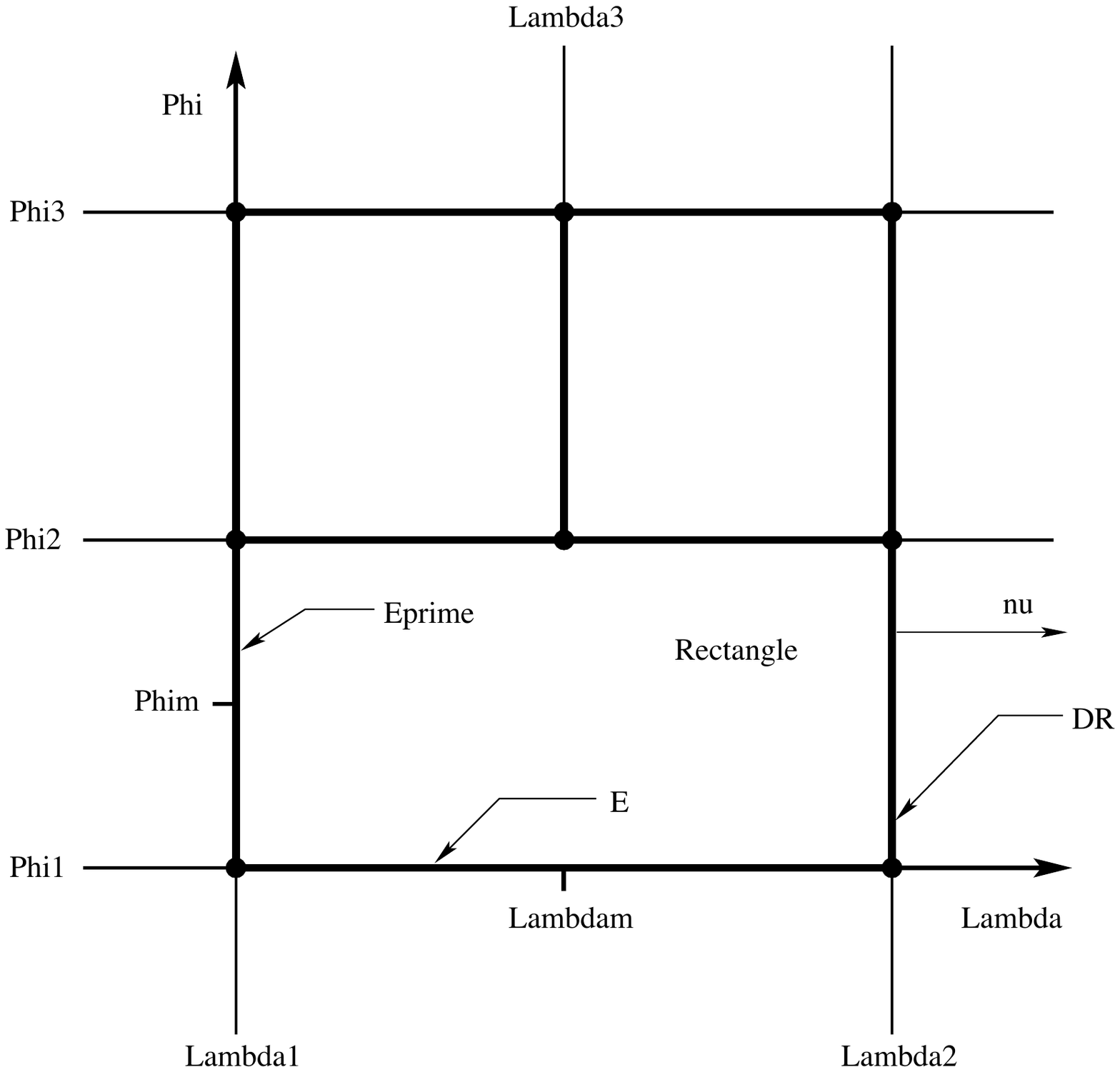}
\caption{Five-sided rectangular cell $\calR$ (on southern
          hemisphere of $\See$)}
\label{fig:cells.2}
\end{center}
\end{figure}

The general structure of our grid is shown in \figref{fig:sphere.1},
and its essential feature is the following.
Every cell $\calR$ is bounded by sides which lie either
along a fixed latitude circle ($\phi=const.$) or a fixed longitude
circle ($\lambda=const.$).
We have
\be
\begin{split}
\calR& := \big\{\lambda_1 \leq\lambda\leq\lambda_2, \quad
\phi_1 \leq\phi\leq\phi_2 \big\},
\label{eq:sphere.12}
\end{split}
\ee
as represented in \figref{fig:cells.1}. In most cases, $\del\calR$
consists of the four sides of $\calR$.
However, across special latitude circles we reduce the number of cells,
so that the situation (for a reduction by ratio of $2$) is as in
\figref{fig:cells.2}.
In this case the boundary $\del\calR$ consists of five sides,
(so that the intermediate point $(\lambda_3,\phi_2)$ is regarded
 as an additional vertex), and even in this five-sided cell $\calR$
every side satisfies the above requirement.

\noindent
The length of a side
$\,e=\big\{\lambda_1\leq\lambda\leq\lambda_2,\,\phi=const.\big\}$
equals $(\lambda_2-\lambda_1)\cos\phi$, while the length of a
side $\,e'=\big\{\phi_1\leq\phi\leq\phi_2,\,\lambda=const.\big\}$
is $\phi_2-\phi_1$.
Consequently, the area $A_{\calR}$ of the cell $\calR$ is
\be
\begin{split}
A_{\calR}&= \int_{\lambda_1}^{\lambda_2}\diff\lambda
              \int_{\phi_1}^{\phi_2}cos{\phi}\,\diff\phi =
             (\lambda_2-\lambda_1)(\sin{\phi_2}-\sin{\phi_1}).  
\end{split}
\nonumber
\ee

\

\

\subsection{Geometry-compatible discretization of the divergence
            operator}
\label{discrete_divergence.sec}
\noindent
Given any rectangular domain $\calR$ of the form \eqref{eq:sphere.12},
the approximate flux divergence is now derived as an approximation
of the integral of the flux
along the boundary $\del\calR$, divided by its area, as follows:
\be
\aligned
& \Bigl( \nabla_{\scscT}\cdot\bfF(\bfx,u) \Bigr)^{approx}
=
   \frac{I_{\calR}}{A_{\calR}},
   \qquad
I_{\calR} = \Bigl( \oint\limits_{\del\calR}\,
              \bfF(\bfx,u)\cdot\bfnu\, ds\Bigr)^{approx},
\endaligned
\label{eq:S2-flux-rectangle}
\ee
where $ds$ is the arc length along $\del\calR$ and $\bfnu$ is
 the outward-pointing unit normal to $\del\calR \subset \Se^2$.
In the limit $\lambda_2,\phi_2\rightarrow\lambda_1,\phi_1$
the approximation \eqref{eq:S2-flux-rectangle} to the divergence
term approaches the exact value \eqref{eq:sphere.5}.

We need to check that the geometric compatibility condition
\eqref{eq:sphere.9} is satisfied for the approximate flux
divergence.
This requirement will be taken into account in formulating
our finite volume scheme for \eqref{eq:sphere.6}. 


Consider now the actual evaluation of the term $I_{\calR}$ defined in
\eqref{eq:S2-flux-rectangle} and consider the cell shown in
\figref{fig:cells.1}, under the assumption that $u=u(\lambda,\phi,t)$
is smooth on ${\calR}$.
We propose to approximate the flux integral along each edge of
${\calR}$ in the following way.
As in Section~\ref{constr}, let us decompose the flux into its
$(\lambda,\phi)$ components:~
$$
\bfF(\bfx,u)=F_\lambda(\lambda,\phi,u)\bfi_\lambda+
              F_\phi(\lambda,\phi,u)\bfi_\phi.$$
On each side the integration is carried out by (i) taking midpoint
values of the appropriate flux component, and (ii) using the correct
arc-length of the side.
We designate the midpoints of the edge $\,e\,$ as
$\,\lambda^{e,m}=(\lambda_1+\lambda_2)/2$ and
$\phi^{e,m}=\phi_{1}$ (see \figref{fig:cells.1}),
and likewise for the edge $\,e'$.

Throughout the rest of this section we {\sl restrict attention
to the gradient flux vector} constructed in Claim~\ref{cor2.2}.
In particular, it comprises the class of {\sl homogeneous flux vectors,}
given by
\eqref{eq:sphere.10}--\eqref{eq:sphere.11}.

Taking $u$ as constant $u=u^{e,m}$ along the
side $e\in \del\calR$, the total approximate flux is given by
\begin{equation}
\label{eq:flux-on-e}
\Big[\;\oint\limits_{e}\,\bfF(\bfx,u)\cdot\bfnu\,\diff s
  \Big]^{approx}=-\big( h(e^{2},u^{e,m})-h(e^{1},u^{e,m}) \big),
\end{equation}
where $e^{1},e^{2}$ are, respectively, the initial and final
endpoints of $e$ (with respect to the sense of the integration).

\noindent
Summing up over all edges we obtain:

\

\begin{Claim}[Discrete geometry-compatibility condition.]
Consider the gradient flux vector constructed in
Claim \ref{cor2.2}. Then, if $u\equiv const.$, $I_{\calR}=0$, so that
$$
\bigl[\nabla_{\scscT}\cdot\bfF(\bfx,u)\bigr]^{approx}=0,
$$
and thus a discrete version of
the divergence-free condition \eqref{eq:sphere.9} holds.
\end{Claim}

\

\begin{Remark} The claim above applies to gradient flux vectors in Claim \ref{cor2.2},
and, in particular, to homogeneous flux
\eqref{eq:sphere.10}--\eqref{eq:sphere.11}.
On the other hand, for a more general geometry-compatible flux
$\bfF(\bfx,u)$, such a result can be obtained
only if the dependence on $\bfx$ is integrated exactly along each side,
a requirement that must be imposed on the scheme.
\end{Remark}

\subsection{Godunov-type approach to the numerical flux}
\label{Evaluation_of_approximate_fluxes.sec}
We continue to deal with the gradient flux given in
Claim \ref{cor2.2}.
We assume different (constant) values of
$u=u(\lambda,\phi,t)$ in grid cells  and evaluate the numerical flux
values at each edge
from the solution to a Riemann problem with data comprising these
values $u(\lambda,\phi,t)$ in the cells on either side of that edge.
At the midpoint
$(\lambda^{e,m},\phi^{e,m})$ of each side $e$
we solve the Riemann problem in a direction perpendicular to $e$,
and denote the resulting solution $u^{e,m}$.
The corresponding fluxes are then evaluated as
$\bfF(\lambda^{e,m},\phi^{e,m},u^{e,m})$.

\begin{figure}[!htb]
\psfrag{Lambda1}[][ml]{\hspace{ 2.8em}  {$\lambda_1$}}
\psfrag{Lambda2}[][ml]{\hspace{ 1.5em}  {$\lambda_2$}}
\psfrag{Lambda3}[][tl]{\hspace{ 2.7em}  {$\lambda_3$}}
\psfrag{Phi1}[][ml]{\hspace{ 1.6em}  {$\phi_1$}}
\psfrag{Phi2}[][ml]{\hspace{ 2.0em}  {$\phi_2$}}
\psfrag{Phim}[][tl]{\hspace{ 2.2em}  {$\phi^{e',m}$}}
\psfrag{u=u_L}[][bl]{\hspace{ 1.0em}  {$u=\uL$}}
\psfrag{u=u_R}[][bl]{\hspace{ 3.0em}  {$u=\uR$}}
\psfrag{u=u^e,m}[][tl]{\hspace{ 2.5em}  {$u=u^{e',m}$}}
\psfrag{M}[][bl]{\hspace{ 0.7em} {M}}
\psfrag{Eprime}[][tl]{\hspace{ 0.3em} {$e'$}}
\begin{center}
\includegraphics[clip,width =106mm,bb=60 325 530 520]
                                              {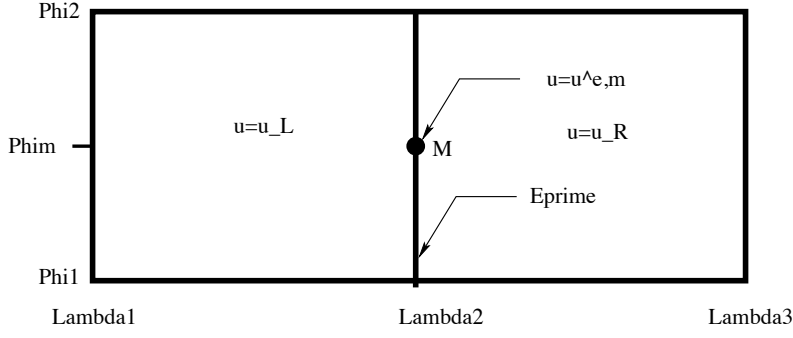}
\caption{Two $\lambda$-adjacent cells with constant states
             $\uL,\,\uR$}
\label{two_cells_lamda.1}
\end{center}
\end{figure}

\begin{figure}[!htb]
\psfrag{Lambda1}[][bl]{\hspace{ 3.3em}  {$\lambda_1$}}
\psfrag{Lambdam}[][bl]{\hspace{ 4.5em}  {$\lambda^{e,m}$}}
\psfrag{Lambda2}[][bl]{\hspace{ 2.2em}  {$\lambda_2$}}
\psfrag{Phi1}[][ml]{\hspace{ 1.6em}  {$\phi_1$}}
\psfrag{Phi2}[][ml]{\hspace{ 2.0em}  {$\phi_2$}}
\psfrag{Phi3}[][ml]{\hspace{ 2.5em}  {$\phi_3$}}
\psfrag{u=u_L}[][ml]{\hspace{ 2.0em}  {$u=\uL$}}
\psfrag{u=u_R}[][ml]{\hspace{ 2.0em}  {$u=\uR$}}
\psfrag{u=u^e,m}[][ml]{\hspace{ 4.3em}  {$u=u^{e,m}$}}
\psfrag{M}[][bl]{\hspace{ 0.2em} {M}}
\psfrag{E}[][ml]{\hspace{ 0.3em} {$e$}}
\begin{center}
\includegraphics[clip,width=60mm,bb=155 242 450 602]
                                                {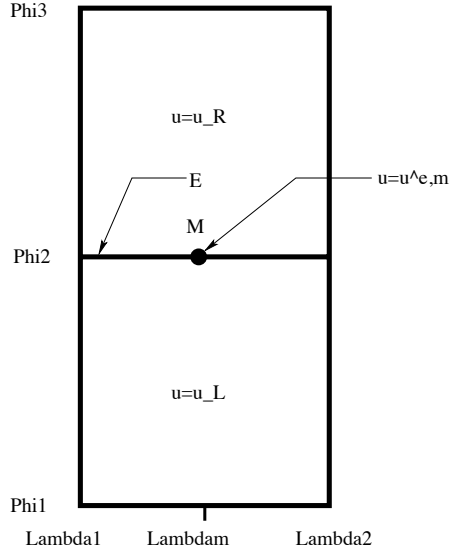}
\caption{Two $\phi$-adjacent cells with constant states
             $\uL,\,\uR$}
\label{two_cells_phi.1}
\end{center}
\end{figure}

We can split Eq.~\eqref{eq:sphere.6}
by invoking the explicit form of the divergence
\eqref{eq:sphere.5},
getting
\be\begin{split}
\label{eq:split-lambda.phi}
  \frac{\del u}{\del t} +
  \frac{1}{\cos{\phi}}
 \frac{\del}{\del\lambda} F_{\lambda}(\lambda,\phi,u) = 0
 \qquad
&\text{for the side} 
\qquad e' : \lambda=\lambda_{2} , \qquad \EQSPLITLAM\nonumber
\\
  \frac{\del u}{\del t} -
  \frac{1}{\cos{\phi}}
 \frac{\del}{\del\phi}\Bigl( F_{\phi}(\lambda,\phi,u)
   \cos{\phi}\Bigr) = 0
    \qquad
&\text{for the side} 
\qquad e\, : \phi=\phi_{2}, \qquad \EQSPLITPHI\nonumber
\end{split}\ee
\stepcounter{equation}
Consider two adjacent cells, as in
\figref{two_cells_lamda.1}
or in
\figref{two_cells_phi.1}.
By fixing $\phi=\phi^{e,m}$ (resp. $\lambda=\lambda^{e,m}$) in
\EQSPLITLAM
(resp. \EQSPLITPHI) we can evaluate $u=u^{e,m}$ as a one-dimensional
solution at $\lambda=\lambda^{e,m}$ (resp. $\phi=\phi^{e,m}$).

\

\

We include here some remarks that will be useful in the implementation of the scheme.

Consider an homogeneous flux vector as in Claim \ref{cor2.1}
so that its components are given by
\eqref{eq:sphere.11}.
Suppose that $u(\lambda,\phi,t_{n})=\uL$
(resp. $u(\lambda,\phi,t_{n})=\uR$) in the cell
$\big\{ \lambda_1<\lambda<\lambda_2,\,\phi_1<\phi<\phi_2\big\}$
(resp. $\big\{ \lambda_2<\lambda<\lambda_3,\,\phi_1<\phi<\phi_2 \Big\}$),
as in
\figref{two_cells_lamda.1}.
At the point
\hbox{$M(\lambda^{e',m},\phi^{e',m})$}
Eq.~\EQSPLITLAM
takes the form
\be
\begin{split}
  \frac{\del u}{\del t}
&+\tan{\phi^{e',m}}
  \frac{\del}{\del\lambda}
\Bigl( f_{1}(u)\cos{\lambda}+ f_{2}(u)\sin{\lambda}\Bigr) -
\frac{\del}{\del\lambda} f_{3}(u) = 0.
\label{eq:SCL.S2-lambda.phi_m}
\end{split}
\ee
Setting
\be
\begin{split}
  g(\lambda,u) &=
  \tan{\phi^{e',m}}\Bigl(
        f_{1}(u)\cos{\lambda}+
        f_{2}(u)\sin{\lambda}       \Bigr) - f_{3}(u),
\label{eq:SCL.S2-g(u)}
\end{split}
\ee
we see that equation~\eqref{eq:SCL.S2-lambda.phi_m}
is the scalar one-dimensional conservation law
\be
\begin{split}
  \frac{\del u}{\del t} +
\frac{\del}{\del\lambda}g(\lambda,u) &= 0, \hspace{33pt}
                                            t\geq t_{n}
\label{eq:SCL.S2-scalar.lambda.g(u)}
\end{split}
\ee
subject to the initial data $u=\uL$ (resp. $u=\uR$)
for $\lambda<\lambda_2$ (resp. $\lambda>\lambda_2$).

Likewise, we repeat the former analysis for $\phi$-adjacent cells
by taking the constant states
$u(\lambda,\phi,t_{n})\EQU\uL$, $u(\lambda,\phi,t_{n})\EQU\uR$
in cells $\big\{ \lambda_1\LT\lambda\LT\lambda_2,\,
        \phi_1\LT\phi<\phi_2\big\}$,
$\big\{ \lambda_1\LT\lambda\LT\lambda_2,\,
        \phi_2\LT\phi\LT\phi_3 \big\}$,
as depicted in
\figref{two_cells_phi.1}.
At the point
\hbox{$M(\lambda=\lambda^{e,m},\,\phi=\phi_2)$},
the equation
\EQSPLITPHI
then takes the form
\be\begin{split}
  \frac{\del u}{\del t} &+ \frac{1}{\cos{\phi}}
  \frac{\del}{\del\phi}\Bigl(
 -\sin{\lambda^{e,m}}\cos{\phi}f_{1}(u) +
  \cos{\lambda^{e,m}}\cos{\phi}f_{2}(u)\Bigr) = 0.
\label{eq:SCL.S2-lambda_m.phi}
\end{split}
\ee
We then set the $\phi$-flux function
\be
\begin{split}
  k(\phi,u) &= \Bigl(
 -\sin{\lambda^{e,m}}\,f_{1}(u) +
  \cos{\lambda^{e,m}}\,f_{2}(u) \Bigr) \cos{\phi}
\label{eq:SCL.S2-h(u)}
\end{split}
\ee
so that
equation~\eqref{eq:SCL.S2-lambda_m.phi}
is the scalar one-dimensional conservation law
\be
\begin{split}
  \frac{\del u}{\del t} + \frac{1}{\cos{\phi}}
\frac{\del}{\del\phi}k(\phi,u) &= 0, \hspace{33pt}
                                            t\geq t_{n}
\label{eq:SCL.S2-scalar.phi.k(u)}
\end{split}
\ee
subject to the initial data $u=\uL$ (resp. $u=\uR$)
for $\phi<\phi_2$ (resp. $\phi>\phi_2$).

\subsection{Solution to the Riemann problem}
\label{Solution.to.the.Riemann.Problem}
The solution at the discontinuity $\lambda\EQU\lambda_2$
at the initial time $t=t_{n}$ is given by the Riemann solution to
\EQSPLITLAM.
For simplicity of the presentation we specialize here to the flux
\eqref{eq:SCL.S2-scalar.lambda.g(u)}.
Since the dependence of $g(\lambda,u)$ on $\lambda$ is smooth, this
solution is obtained by fixing $\lambda=\lambda_{2}$, thus solving the
classical conservation law
\be
\begin{split}
  \frac{\del u}{\del t} +
\frac{\del}{\del\lambda}g(\lambda_{2},u) &= 0, \hspace{33pt}
                                            t\geq t_{n}
\label{eq:SCL.S2-scalar.lambda.g(u).2}
\end{split}
\ee
subject to the initial jump discontinuity of $u$.

We denote this solution by $u^{2,m}$.
Observe that the flux $g(\lambda,u)$ in
\eqref{eq:SCL.S2-scalar.lambda.g(u).2}
is in general non-convex.
The Riemann solution may therefore consist of several waves.
It is a self-similar solution depending only on the slope
$(\lambda-\lambda_2)/(t-t_{n})$.
The value $u^{2,m}$ is the value along the line
$\lambda\EQU\lambda_2$.
It therefore corresponds either to a sonic wave, namely
$g'(\lambda_{2},u^{2,m})\EQU 0$, or to an ``upwind value''
$u\EQU\uL$ (resp. $u\EQU\uR$)
in the case where all waves propagate to the right (resp. left).

Actually, the procedure for solving the Riemann problem in the case
of a nonconvex flux function $g(\lambda_{2},u)$ is well-known and goes
back to classical works by Oleinik and others.
We recall it here briefly.  Assume first that $\uL\LT\uR$.
Consider the {\sl convex envelope} of $g$,
namely, the largest convex continuous  function $g_c$, over the
interval $[\uL,\uR]$, such that $g_c\LEQ g$ at all points.
Clearly, $g_c\EQU g$ in ``convex sections'' of the graph
of $g$, while it consists of linear segments when $g_c\LT g$.
It is easy to see that
the ``convex segments'', where $g\EQU g_{c}$,
represent rarefaction waves
(in the full Riemann solution) while the linear segments represent
jumps (i.e., shock waves).
In particular, the solution $u^{2,m}$ is given by the
following formula:
\be
u^{2,m}=v_{min},\quad \text{where}\quad
        g(\lambda_{2},v_{min})\leq g(\lambda_{2},v) \quad
        \text{for all}\quad  v\in[\uL,\uR].
\label{eq:SCL.S2-nonconvex-RP.u2m}
\ee
There are in fact three possibilities for this solution:
\begin{description}
\item[a)]
  $\uL<u^{2,m}<\uR$, which implies that
  $g'(\lambda_{2},u^{2,m})=0$ (a sonic point).
\item[b)]
  $u^{2,m}=\uL$, the whole wave pattern moves to the right.
\item[c)]
  $u^{2,m}=\uR$, the whole wave pattern moves to the left.
\end{description}

\

Similarly, in the case $\uL\GT \uR$, we construct the
``concave envelope'' of $g$, namely,
the smallest concave continuous function $g_{c}$ such that
$g_{c}\GEQ g$.
Again the linear segments correspond to jump discontinuities while the
concave segments ($g\EQU g_{c}$) correspond to rarefaction waves.
The solution to the Riemann problem is now given by
$u^{2,m}\EQU v_{max}$, where
$g(\lambda_{2},v_{max})\GEQ g(\lambda_{2},v)$, $v\in [\uR,\uL]$.
As above, there are three possibilities for the solution
(sonic, left-upwind, or right-upwind).

Replacing in the foregoing analysis the $\lambda$-flux function
$g(\lambda_{2},u)$ by the $\phi$-flux function $k(\phi_{2},u)$,
the equation \eqref{eq:SCL.S2-scalar.phi.k(u)} reads
\be
\begin{split}
  \frac{\del u}{\del t} + \frac{\del}{\del\phi}\Bigl(
       -\sin{\lambda^{2,m}}f_{1}(u) + \cos{\lambda^{2,m}}f_{2}(u)\Bigr)
                            &= 0, \hspace{33pt} t\geq t_{n}\,.
\label{eq:SCL.S2-scalar.phi.k(u).2}
\end{split}
\ee
We get the Riemann solution to
\eqref{eq:SCL.S2-scalar.phi.k(u).2}
in the three cases {\bf a)}, {\bf b)}, {\bf c)} as above.

\subsection{Convergence proof}
\label{Proof.of.convergence.of.the.scheme}
The computational elements (``grid cells'') are denoted in~\cite{ABL}
by $K$.
Their sides are denoted by $e$ and the flux
function across $e$ is given by $f_{e,K}(u,v),$ where $u$
is the (constant) value in $K$ and $v$ is the value in the
neighboring cell (sharing the same side $e$) $K_e.$
In our grid of the sphere, some cells are actually {\sl pentagons;}
these are the cells whose lower-latitude side (along a latitude
$\phi=const$) borders the two higher-latitude sides of the two
lower-latitude neighbor cells, as shown in
\figref{fig:cells.2} for the southern hemisphere grid.
For such cells, the lower-latitude
side consists of {\sl two} faces, each one of them common with
one of the lower-latitude neighboring cells.

With this construction of the grid, we can check the conditions in~\cite{ABL}
imposed on the numerical flux. It is important to keep in mind that we are
dealing with the gradient flux vectors given by Claim \ref{cor2.2}.

\

\begin{Claim}[Convergence of the proposed scheme.]
Consider the first-order finite volume scheme described above.
Assume that the flux vector has the gradient form in Claim \ref{cor2.2}. Let
 $f_{e,\calR}(u,v)$ be the numerical flux calculated on the side
$e$ of the computational cell
 $\calR,$ using
\eqref{eq:flux-on-e}, where the midpoint value of $u$ is obtained from
the Riemann solution.
Then $f_{e,\calR}(u,v)$ satisfies the assumptions
(5.5)-(5.7)
of \cite{ABL},
and the numerical solution converges to the exact solution as the
maximal size of the grid cells shrinks to zero.
\end{Claim}

\

\noindent{\sl Proof.}
Consider the flux across a longitude side $e:\lambda=\lambda_2$,
which is given by $F_\lambda$ in the equation~\EQSPLITLAM.
The procedure for integrating the flux across $e$
is described by \eqref{eq:flux-on-e},
while in Subsection
\ref{Solution.to.the.Riemann.Problem} the calculation of
$F_{\lambda}(\lambda_{2},\phi^{2,m},u^{2,m})$ is described.
It can be summarized as follows.

First, the solution $u^{2,m}$ to the Riemann problem associated with
equation~\EQSPLITLAM is found, assuming $u,v$ to be the values on the
two sides.
However, note that $F_\lambda$ depends {\sl explicitly} on $\phi$,
and to be precise we need to replace in \EQSPLITLAM the mean value
$\phi^{2,m}$ by $\phi.$
Thus, we find $u^{2,m}=u^{2,m}(\phi).$

Clearly, in the case $u=v$ we get identically
$u^{2,m}(\phi)=u=v$ and so the exact flux satisfies
$$F_\lambda=F_\lambda(\lambda_2,\phi,u^{2,m})$$
and its integration will give exactly the approximate value
$$
f_{e,K}(u,v)=-\big( h(e^{2},u^{2,m})-h(e^{1},u^{2,m})\big),
$$
as in \eqref{eq:flux-on-e}.
Thus, condition (5.5) in~\cite{ABL} is satisfied.

Clearly, the conservation property (5.6) is satisfied even
with the approximate definition.

Also, the flux as defined in \eqref{eq:flux-on-e} makes it easy to
check (5.7), as the flux is independent of $\phi$ and the
monotonicity is thus a result of general
properties of the Riemann solver (even for nonconvex fluxes).
For example, if $u<v$, one considers the convex envelope of
$F_\lambda$, as defined in \EQSPLITLAM
(with $\phi=\phi^{2,m}$) and then considers $u^{2,m}$ as the
{\sl minimal} value on this envelope (over $[u,v]$).
Clearly changing $u$ upward will either change $u^{2,m}$
upward or leave it unchanged.
This completes the proof.
\qed

%
\section{Second-order extension based on the GRP solver}
\label{second}
To improve the order of accuracy, we consider again the cell
$\;\lambda_1\LT \lambda\LT \lambda_2,\,
      \phi_1\LT    \phi\LT    \phi_2\;$
and assume that $u$ is \emph{linearly} distributed there.
We use $\uLlambda,\uLphi$
(resp. $\uRlambda,\uRphi$) to denote the slopes in the cell to
the left (resp. right) of the  side $\lambda\EQU\lambda_2$.
We also denote by $\uL(\phi)$ (resp. $\uR(\phi)$) the limiting value
 (linearly distributed) of $u$ at $\lambda\EQU\lambda_2-$
(resp. $\lambda\EQU\lambda_2+$).
Clearly, the solution to the Riemann problem
across the discontinuity is a function of $\phi$, and we denote it
by $u^{2,m}(\phi)$, which conforms to our notation in
Subsecion~\ref{Solution.to.the.Riemann.Problem} above
(where $u$ was constant on either side of the discontinuity).
The value of $u^{2,m}(\phi)$ is obtained by solving the Riemann
problem associated with
Eq.~\EQSPLITLAM
with $\phi^{2,m}$ replaced by $\phi$, subject to the initial data
$\uL(\phi),\; \uR(\phi)$.
Restricting to the middle point
$\phi=\phi^{2,m}$, the solution $u^{2,m}(\phi^{2,m})$
(at $\lambda=\lambda^{2,m}$) is in one of
the three categories listed above (i.e., sonic, left-upwind, right-upwind).
By continuity, the solution
$u^{2,m}(\phi)$ will still be in the same category for $\phi-\phi^{2,m}$
sufficiently small.
The solution at $(\lambda^{2,m},\phi^{2,m})$ varies in time
and the GRP method
deals with the determination of its time-derivative at that point.

Accounting for the variation of the solution over a time interval enables us to modify
 the Godunov approach to the determination of edge fluxes , as presented in
Section~\ref{Evaluation_of_approximate_fluxes.sec}. We assume that  the flux vector
depends explicitly  on $\bfx,$ as in \eqref{eq:sphere.8}. In what follows we use for
simplicity the ``imbedded'' notation $\bfx=(x_1,x_2,x_3)$ for a point on the sphere (see
the Introduction), along with the corresponding spherical coordinates $\lambda, \phi.$
 We further assume that the vector field $\bfPhi$ is given by the following
extension of \eqref{eq:sphere.10}
\be
\begin{split}
\bfPhi(\bfx,u) &= \nabla_{\!\!\bfx} h(\bfx,u)
\\ &=
              q_{1}(x_{1})f_1(u)\,\bfi_1 +
              q_{2}(x_{2})f_2(u)\,\bfi_2 +
              q_{3}(x_{3})f_3(u)\,\bfi_3,
\label{eq:h-scheme.3}
\end{split}
\ee
The zero-divergence identity is  obtained as a result of expressing
$\bfPhi$ as a gradient $\nabla h$ in the sense of Claim \ref{cor2.2}.

For our choice of $\bfPhi$ such a  representation of $\bfPhi$ as gradient of $h$ is
obtained when $h$ is taken as
\be
\begin{split}
h(\bfx,u) &=  r_{1}(x_{1})f_1(u) +
              r_{2}(x_{2})f_2(u) +
              r_{3}(x_{3})f_3(u)\,,
\label{eq:h-scheme.4}
\end{split}
\ee
and $\,q_{j}(x_{j})=r_{j}'(x_{j}),\,\; j=1,2,3.\,$

Using \eqref{eq:sphere.5} together with the geometry-compatibility property, we get
an explicit form of the conservation law ~\eqref{eq:sphere.6} in our case as
  \be\begin{split}
  \frac{\del u}{\del t} &-
  \sin{\lambda}q_{1}(x_{1})\frac{\del}{\del\phi}f_1(u) +
  \cos{\lambda}q_{2}(x_{2})\frac{\del}{\del\phi}f_2(u) \\
&+\tan{\phi}\Big(
  \cos{\lambda}q_{1}(x_{1})\frac{\del}{\del\lambda}f_1(u) +
  \sin{\lambda}q_{2}(x_{2})\frac{\del}{\del\lambda}f_2(u)\Big)
   -q_{3}(x_{3})
 \frac{\del}{\del\lambda} f_3(u) = 0.
\label{eq:h-scheme.2}
\end{split}
\ee
 The numerical approximation to this equation requires an {\sl operator splitting approach,}
where the derivatives with respect to  $\phi$ and $\lambda$ are considered separately.
We note that such a splitting has already been implemented in the Godunov case,
~\eqref{eq:split-lambda.phi}, in the most general case. In that case, no use has been
made of the geometry-compatibility property. Indeed, this has no bearing on the
\textit{first-order} scheme since the solution to the Riemann problem is obtained by
``freezing" the explicit dependence on $\lambda, \phi$ (and, in particular, ignoring
the terms involving the derivatives with respect to this explicit dependence).\newline
In the present (second-order) situation we proceed as follows.

The ``$\lambda$-split'' equation obtained from \eqref{eq:h-scheme.2}, is
\be
\begin{split}
  \frac{\del u}{\del t}
&+\tan{\phi^{2,m}} \Big( q_{1}(x_{1})\cos{\lambda}
   \frac{\del}{\del\lambda}f_1(u)+
q_{2}(x_{2})\sin{\lambda}\frac{\del}{\del\lambda}f_2(u)\Big) -
q_{3}(x_{3})\frac{\del}{\del\lambda} f_3(u) = 0. \label{eq:h-scheme.5}
\end{split}
\ee
Note that the coefficients are retained as functions of $\lambda$ and
are not ``frozen'' at $\lambda=\lambda^{2,m}$.
This is of course due to the fact that in employing
the GRP scheme we consider $\lambda$-derivatives on
either side of the edge, so as in any limiting analysis, we must first let
$\lambda \to \lambda^{2,m}$, then substitute $\lambda=\lambda^{2,m}$.

The $\lambda$-edge flux function $g(\lambda,u)$ (compare \eqref{eq:SCL.S2-g(u)}), is now
extended to $g(\bfx,u)$ as
\be
\begin{split}
  g(\bfx,u) &=
  \tan{\phi^{2,m}}\Bigl(
        q_{1}(x_{1})\cos{\lambda}f_1(u)+
        q_{2}(x_{2})\sin{\lambda}f_2(u)  \Bigr) - q_{3}(x_{3})f_3(u),
\end{split}
\nonumber
\ee
and the scalar one-dimensional conservation law under consideration is now
rewritten as an equation with a source term (a balance law) 
\be
\begin{split}
& \frac{\del u}{\del t} + \frac{\del}{\del\lambda}g(\bfx,u) =
         S_{\lambda}\,,
         \hspace{33pt} t > t_{n} \\
& S_{\lambda} = \tan{\phi^{2,m}}\Bigl(f_1(u)\frac{\del}{\del\lambda}\big(
      q_{1}(x_{1})\cos{\lambda} \big) +
                    f_2(u)\frac{\del}{\del\lambda} \big(
      q_{2}(x_{2})\sin{\lambda} \big) \Bigr) -
      f_3(u)\frac{\del}{\del\lambda}q_{3}(x_{3}),
\label{eq:h-scheme.7}
\end{split}
\ee
subject to the initial data (for $u$ and its slope) $\uL(\phi^{2,m})$,
$\uLlambda$ (resp.~$\uR(\phi^{2,m})$, $\uRlambda$)
for $\lambda<\lambda_2$ (resp.~$\lambda>\lambda_2$).
Observe that the equation is written in a ``quasi-conservative form", which offers more
convenience in the GRP treatment
\cite[Chap.~5]{grp_book}.
The right-hand side term $S_{\lambda}$ is
just the result of the $\lambda$ differentiation of the flux $g(\bfx,u).$ Obviously,
the geometry-compatibility condition implies that this source term should cancel out with the
corresponding source term in the ``$\phi$-split" equation.
The solution $u^{2,m}$ to the Riemann
problem is obtained by freezing the coordinate $\lambda$ at its
edge value, so that, in particular,  the source term in
\eqref{eq:h-scheme.7} can be taken as zero at this stage.

In the framework of the GRP analysis, the source term $S_{\lambda}$ is added to terms
arising from the piecewise-linear initial data, in producing the time-derivative of
the solution $u^{2,m}(\phi^{2,m}) + \frac{\del u}{\del
t}(\lambda^{2,m},\phi^{2,m},t_{n}+) \frac{\Delta t}{2},\;\Delta t = t_{n+1}-t_{n}$. As
explained above, $u^{2,m}(\phi^{2,m})$,  the solution to the associated Riemann problem, is
obtained by using the ``edge values" $\uL(\phi),\; \uR(\phi)$. It remains, therefore,
to determine the instantaneous time-derivative $\frac{\del u}{\del
t}(\lambda^{2,m},\phi^{2,m},t_{n}+)$, as is outlined below.

The time-derivative of $u$  is  given by \be
\begin{split}
\frac{\del u}{\del t}(\lambda^{2,m},\phi^{2,m},t_{n}+)= -u_{m,\lambda}\,
\frac{\del}{\del u}g(\bfx,u)|_{\lambda^{2,m},\phi^{2,m},u^{2,m}},
\end{split}
\nonumber \ee where the slope value  $u_{m,\lambda}$ is obtained by ``upwinding'',
determined by the associated Riemann problem as follows (we start with the ``easy''
categories {\bf b)}, {\bf c)} above). 

\

\begin{itemize}

\item[{\bf b)}] $u^{2,m}=\uL(\phi^{2,m})$. Then, the wave moves to the right and we
set 
\be
\begin{split}
u_{m,\lambda}=\uLlambda. 
\end{split}
\nonumber
\ee


\item[{\bf c)}] $u^{2,m}=\uR(\phi^{2,m})$. Then, the wave moves to the left and we set
\be
\begin{split}
u_{m,\lambda}=\uRlambda. 
\end{split}
\nonumber \ee Finally, the first category deals with the sonic case. As noted above,
it remains sonic in the neighborhood of $\phi^{2,m}$, so that we have there
$\frac{\del}{\del u}g(\bfx,u)|_{\lambda^{2,m},\phi^{2,m},u^{2,m}}$. The
time-derivative of $u$ reduces therefore to \be
\begin{split}
\frac{\del }{\del t}u(\lambda_2,\phi^{2,m},t\EQU t_{n}+)=
 0. 
\end{split}
\nonumber
\ee

\end{itemize}


Finally, the ``$\phi$-split'' equation obtained from
\eqref{eq:h-scheme.2}, is treated in analogy with the
``$\lambda$-split'' procedure outlined above.

\section{Numerical tests}
\label{numeri}

\subsection{First test case: equatorial periodic solutions}

\begin{figure}[!htb] 
\psfrag{X}[][tl]{\hspace{1.0em}  {$x = \lambda/2\pi$}}
\psfrag{U(X,T)}[][bl]{\hspace{ 2.5em}  {$u$}}
\begin{center}
\includegraphics[clip,width =80mm,bb=75 58 540 392]
                {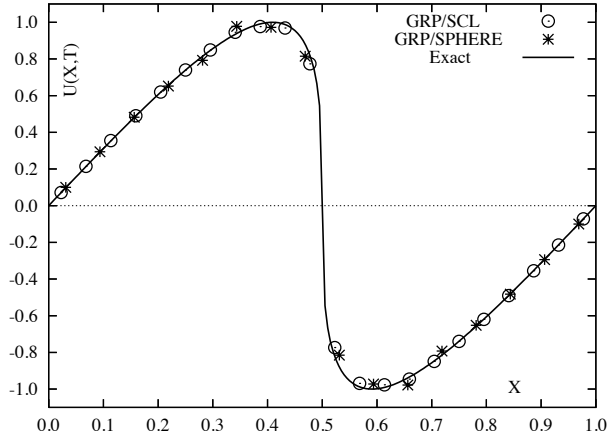}
\caption[Exact, GRP/SCL and GRP/SPHERE solutions to Burgers IVP]
        {Exact, GRP/SCL and GRP/SPHERE ($\Delta\lambda\EQU 2\pi/16$)
         solutions
         to the IVP~\eqref{eq:Burgers-IVP-f_3} at $t=1/2\pi$}
\label{Burgers.1}
\end{center}
\end{figure}

\begin{figure}[!htb]
\psfrag{X}[][tl]{\hspace{1.0em}  {$x = \lambda/2\pi$}}
\psfrag{U(X,T)}[][bl]{\hspace{ 2.5em}  {$u$}}
\begin{center}
\includegraphics[clip,width =80mm,bb=75 58 540 392]
                {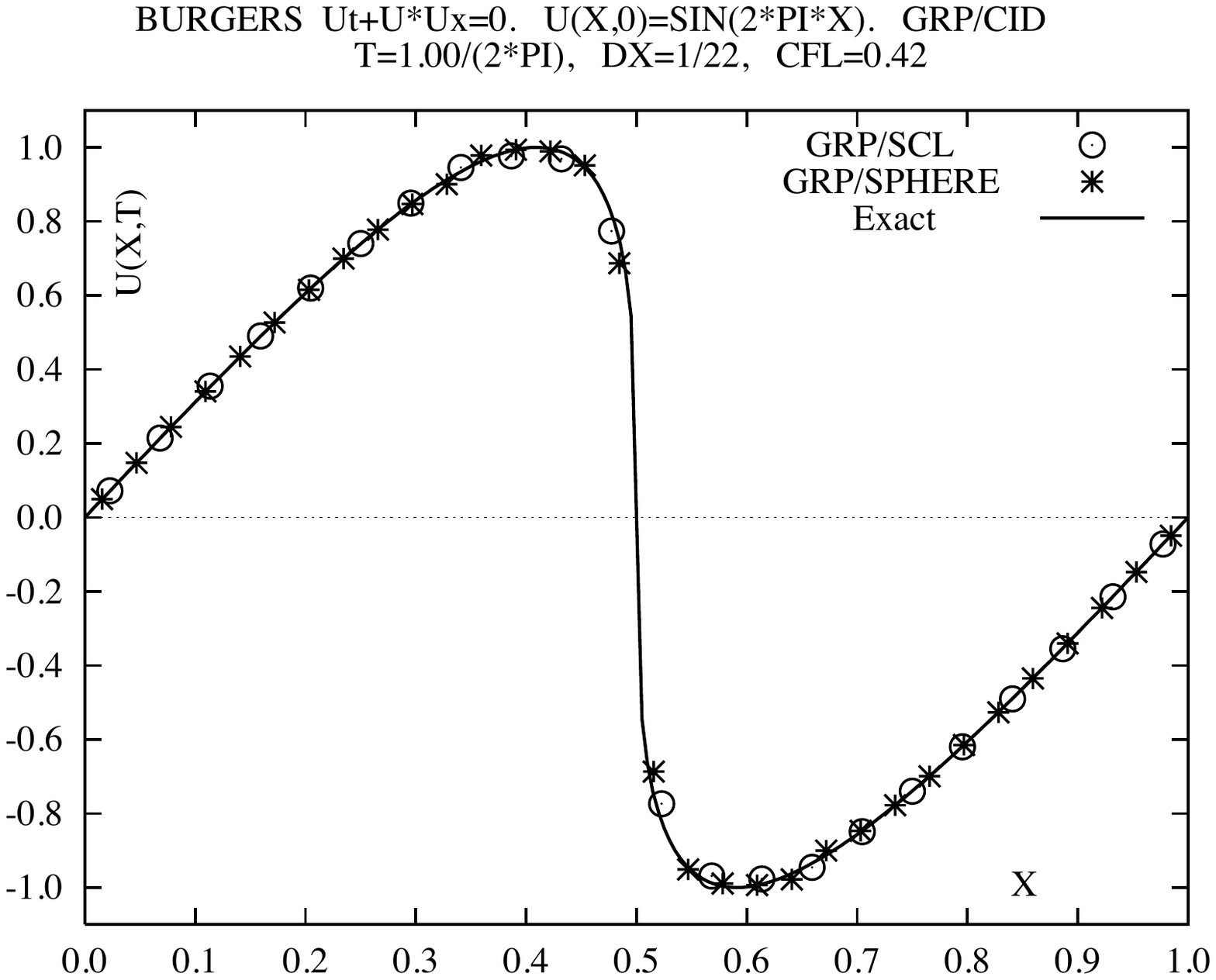}
\caption[Exact, GRP/SCL and GRP/SPHERE solutions to Burgers IVP]
        {Exact, GRP/SCL and GRP/SPHERE ($\Delta\lambda\EQU 2\pi/32$)
         solutions
         to the IVP~\eqref{eq:Burgers-IVP-f_3} at $t=1/2\pi$}
\label{Burgers.2}
\end{center}
\end{figure}

\begin{figure}[!htb]
\psfrag{X}[][tl]{\hspace{1.0em}  {$x = \lambda/2\pi$}}
\psfrag{U(X,T)}[][bl]{\hspace{ 2.5em}  {$u$}}
\begin{center}
\includegraphics[clip,width =80mm,bb=75 58 540 392]
                {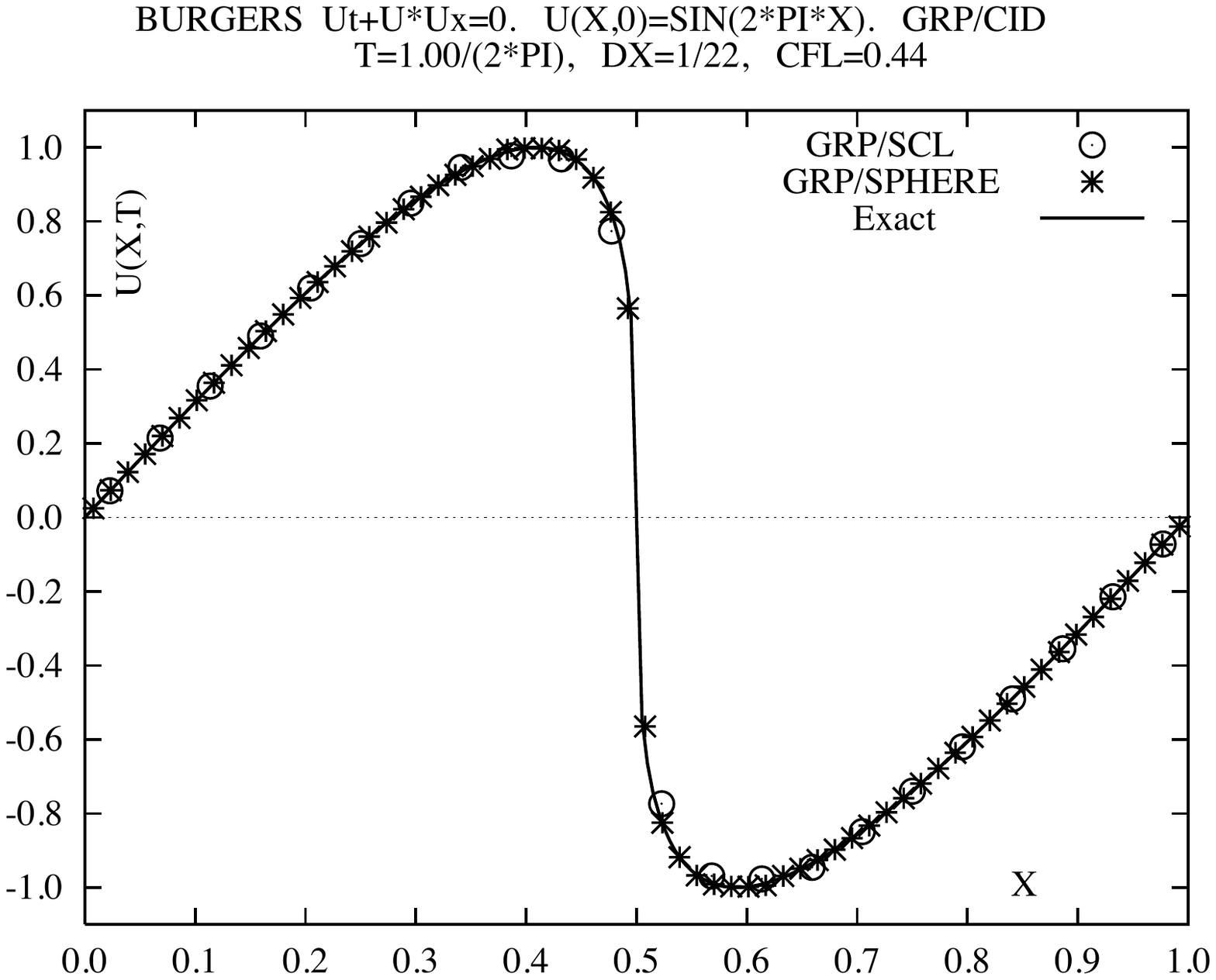}
\caption[Exact, GRP/SCL and GRP/SPHERE solutions to Burgers IVP]
        {Exact, GRP/SCL and GRP/SPHERE ($\Delta\lambda\EQU 2\pi/64$)
         solutions
         to the IVP~\eqref{eq:Burgers-IVP-f_3} at $t=1/2\pi$}
\label{Burgers.3}
\end{center}
\end{figure}

Here, the
conservation law takes the form ~\eqref{equator} and the flux
function and initial data are given by
\be
\begin{split}
f_1(u)&=f_2(u)=0,\quad\qquad f_3(u)=-2\pi\,(u^{2}/2),  \\
u(\lambda,\phi,0)&=
\begin{cases}
\sin{\lambda}, \qquad \qquad 0<\lambda<2\pi,\,0<\phi<\pi/12,\\
            \hspace{16pt}0, \qquad \qquad \text{otherwise}.
\end{cases}
\label{eq:Burgers-IVP-f_3}
\end{split}
\ee
As discussed in Section~3 (see the discussion of solutions to
~\eqref{equator}) it is clear that
the solution here (as a function of $\lambda$)
is identical to the periodic solution for the Burgers equation in
$\Rone,$ with periodic boundary conditions on $[0<x<2\pi]$.
However, we compute the numerical
solution here on our \textit{spherical} grid, and we need to check
not only that it
conforms with the one-dimensional case but that it does not ``leak''
beyond the band supporting the initial data.
The results {\sl at the shock formation time} $t_{s}=1/2\pi$
are shown in \figref{Burgers.1} for $\Delta\lambda=2\pi/16$,
in \figref{Burgers.2} for
$\Delta\lambda=2\pi/32$ and in \figref{Burgers.3} for
$\Delta\lambda=2\pi/64$.
These GRP solutions to \eqref{eq:SCL.S2-scalar.lambda.g(u)} clearly
converge to the
exact solution with refinement of the $\lambda$ grid, and are
comparable  to the
corresponding solution to the scalar conservation law in
$\Rone$ with $\Delta x=2\pi/22$.

\

\subsection{Second test case: steady state solutions}

We refer to
Corollary ~\ref{1Dsteady} and using the notation there we take the
flux vector and initial data as:
\be\begin{split}
f_1(u)&=u^{2}/2,  \qquad
f_2(u)=f_3(u)=0, \\
u(\lambda,\phi,0)&=
\cos{\lambda}\,\cos{\phi}.
\label{eq:TC.2}
\end{split}\ee
Using the terminology of Corollary ~\ref{1Dsteady} we see that the
initial function is
the ``simplest'' possible function, corresponding to $g(x_1)=x_1.$

\begin{figure}[!htb]
\begin{center}
\includegraphics[clip,width =90mm,bb=35 142 548 657]
                                              {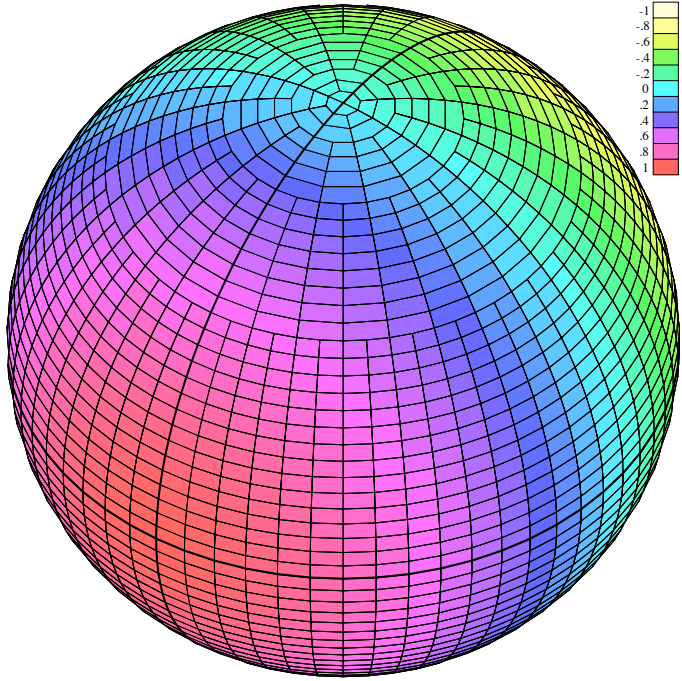}
\caption{Steady-state initial data (and solution)
         to the IVP~\eqref{eq:TC.2} at $t=5$.
\hfill\break\hspace*{28pt}
      Color map range scaled to $\,(u_{min},u_{max})=(-0.998,0.998)$.}
\label{fig:TC.2}
\end{center}
\end{figure}

\

As is shown in \figref{fig:TC.2}, the numerical solution remains
nearly unchanged in
time after being subjected to integration up to $t=5$ by the
GRP scheme with constant time step $\Delta t=0.05$,
the color maps of $u(\lambda,\phi,t)$
at the initial and final times are virtually indistinguishable.
The shown grid has latitude step $\,\Delta\phi=\pi/60,\,$ and an
equatorial longitude step $\,\Delta\lambda=\pi/128$.
A measure $u_{diff}$ to the numerical solution error is defined as the
area-weighted difference $|u(\lambda,\phi,5)-u(\lambda,\phi,0)|$,
obtained by summation over all grid cells.
In this case we obtained $u_{diff}=0.0093$, which is small relative to
the full range $u_{max}-u_{min}=2$.
Hence, the GRP scheme produces an approximation to the
steady-state solution $u(\lambda,\phi,t)=u(\lambda,\phi,0)$ over $\See$.
This test case demonstrates that
the scheme computes correctly the time-evolution for the non-constant
data
\eqref{eq:TC.2}, by calculating an approximately zero
value for the flux divergence in computational cells.


\subsection{Third test case: confined solutions}

We take (as in Claim~\ref{cor2.2})  $\bfPhi(\bfx,u)=\nabla h(\bfx,u)$, where
$h(\bfx,u)=\psi(x_1)x_1 f_1(u).$ The function $\psi(x_1)$ is defined by
\begin{equation}
\psi(x_1)=\begin{cases} 1,\quad x_1\leq 0,\\
              1-6x_1^2+\frac{8}{\sqrt{2}}x_1^3,\quad 0\leq x_1\leq
              \frac{\sqrt{2}}{2},\\
              0,\quad \frac{\sqrt{2}}{2}\leq x_1\,.
\end{cases}
\label{eq:TC.3}
\end{equation}
The flux vector is then given by
 $$\bfF(\bfx,u)=\bfn(\bfx)\times\bfPhi(\bfx,u).$$
The solution is clearly confined to the  sector
$x_1\leq \frac{\sqrt{2}}{2}$ of the sphere.
Its boundary is a circle which intersects the meridian
$\lambda=0$ at $\phi=\frac{\pi}{4}.$

The flux in the subdomain $x_1\leq 0$ is given by
$$\bfF(\bfx,u)=\bfn(\bfx)\times f_1(u)\,\bfi_1,$$
so if we take the initial data as $\psi(x_1)u_0(x_1),$
where $u_0$ is the steady state solution of the
second test case (and also the same $f_1(u)$),
the solution \textsl{remains steady in that part},
namely, in $x_1\leq 0.$ Clearly, it evolves in time in the
region $0\leq x_1\leq \frac{\sqrt{2}}{2}$, but vanishes identically
(for all time) if $\frac{\sqrt{2}}{2}\leq x_1.$

\begin{figure}[!htb]
\begin{center}
\includegraphics[clip,width =90mm,bb=35 142 548 657]
                                              {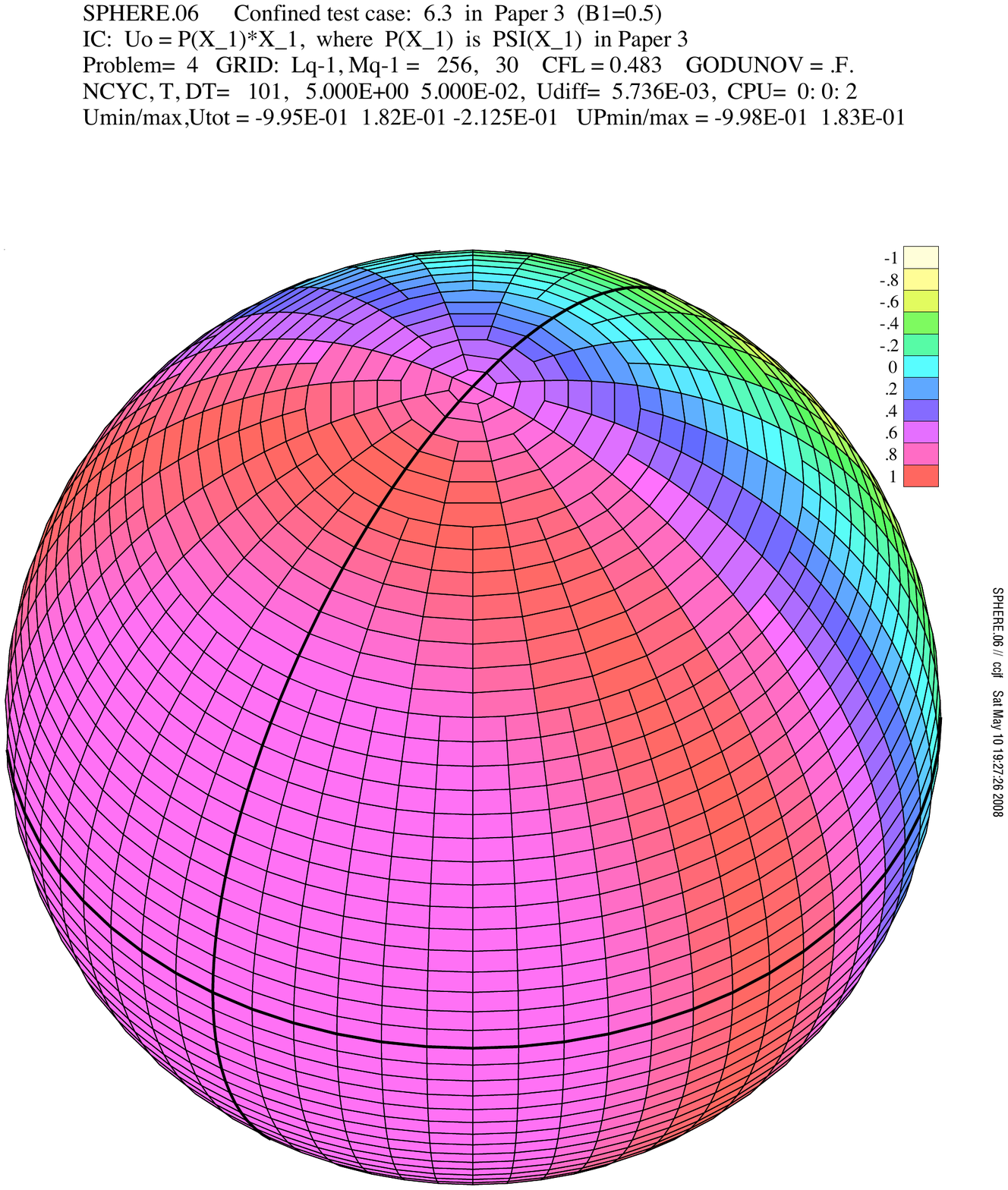}
\caption{Confined solution test case, with the IVP data
         to the IVP~\eqref{eq:TC.3} at $t=5$.
\hfill\break\hspace*{28pt}
      Color map range scaled to $\,(u_{min},u_{max})=(-0.998,0.183)$.}
\label{fig:TC.3}
\end{center}
\end{figure}

The confined IVP was integrated in time up to $t=5$ by the
GRP scheme, using the same grid and time step as in the second test
case (Subsection 6.2).
The solution is represented by the color map in \figref{fig:TC.3}.
Comparing it to the corresponding initial map (not shown here),
it seems nearly unchanged.
In fact, the initial-to-final difference measure obtained is
$u_{diff}=0.0057$, which indicates a nearly steady
solution in the strip $0 < x_{1} < \sqrt{1/2\,}$.
This test case demonstrates that
the scheme computes correctly the time-evolution for the non-constant
``confined'' data \eqref{eq:TC.3}.

\

\section*{Acknowledgments}

The authors were supported by a research grant of cooperation
in mathematics, sponsored by the High Council for Scientific and
Technological Cooperation between France and Israel,
entitled: {\em ``Theoretical and numerical study of geophysical fluid dynamics in general geometry''.}
This research was also partially supported by the A.N.R.
(Agence Nationale de la Recherche) through the grant 06-2-134423 and by the Centre National de la
Recherche Scientifique (CNRS).


\end{document}